\documentclass{amsart}

\theoremstyle{definition}

\theoremstyle{remark}

\numberwithin{equation}{section}

\begin{document}

\title{A product involving the $\beta$-family in stable homotopy theory}

\author{Xiugui Liu}
\address{School of Mathematical Sciences and LPMC,
Nankai University, Tianjin 300071, P. R. China}
\email{xgliu@nankai.edu.cn}
\thanks{The first author was partially supported by the National
Natural Science Foundation of China (Nos. 10501045, 10771105) and
the Fund of the Personnel Division of Nankai University.}

\author{Wending Li}
\address{School of Mathematical Sciences,
Nankai University, Tianjin 300071, P. R. China}
\email{liwending@gmail.com}


\date{November 1, 2009}



\begin{abstract}
In the stable homotopy groups
 $\pi_{q(p^n+p^m+1)-3}(S)$ of the sphere spectrum $S$ localized at the
prime $p$ greater than three, J. Lin constructed an essential family
$\xi_{m,n}$ for $n \geq m + 2 >5$. In this paper, the authors show
that the composite $\xi_{m,n}\beta_{s}\in
\pi_{q(p^n+p^m+sp+s)-5}(S)$ for $2 \leq s < p$ is  non-trivial,
where $q=2(p-1)$ and $\beta_s \in \pi_{q(sp+s-1)-2}(S)$ is the known
$\beta$-family. We show our result by explicit combinatorial
analysis of the (modified) May  spectral sequence.
\end{abstract}

\maketitle

\section{Introduction and statement of results}
Let $A$ be the mod $p$ Steenrod algebra and $S$ the sphere spectrum
localized at an odd prime $p$. To determine the stable homotopy
groups of spheres $\pi_{\ast}(S)$ is one of the central problems in
homotopy theory.

So far, several methods have been found to determine the stable
homotopy groups of spheres. For example we have the classical Adams
spectral sequence (ASS) (cf.~\cite{Adams}) based on the
Eilenberg-MacLane spectrum $K\mathbb{Z}_p$, whose $E_2$-term is
${\rm Ext}_A^{s, t}(\mathbb{Z}_p, \mathbb{Z}_p)$
 and the Adams differential is given by $  d_r: E_r^{s,  t}\rightarrow
 E_r^{s+r,  t+r-1}$.
 We also have the Adams-Novikov spectral sequence
 based on the Brown-Peterson spectrum
 $BP$(cf.~\cite{Miller, Ravenel, Lee}).

Throughout this paper, we fix the prime $p \geq 5$ and set
$q=2(p-1)$. From ~\cite{Liulevicius}, we know that ${\rm
Ext}_A^{1,1}(\mathbb{Z}_p, \mathbb{Z}_p)$ has $\mathbb{Z}_p$-basis
consisting of $a_0\in {\rm Ext}_A^{1,1}(\mathbb{Z}_p,
\mathbb{Z}_p)$, $h_i\in {\rm Ext}_A^{1,p^iq}(\mathbb{Z}_p,
\mathbb{Z}_p)$ for all $i\geq
 0$ and ${\rm Ext}_A^{2,\ast}(\mathbb{Z}_p, \mathbb{Z}_p)$ has
$\mathbb{Z}_p$-basis consisting of $\alpha_2$, $a_0^2$, $a_0h_i$
$(i>0)$, $g_i$ $(i\geq 0)$, $k_i$ $(i\geq 0)$, $b_i$ $(i\geq 0)$,
and $h_ih_j$ $(j\geq i+2,i\geq 0)$ whose internal degrees are
$2q+1$, $2$, $p^iq+1$, $q(p^{i+1}+2p^i)$, $q(2p^{i+1}+p^i)$,
$p^{i+1}q$ and $q(p^i+p^j)$ respectively.

Let $M$ denote the Moore spectrum modulo the prime $p$ given by the
cofibration
\begin{equation} \label{1.1}
S\stackrel{p}{\rightarrow}S\stackrel{i}{\rightarrow}M\stackrel{j}{\rightarrow}\Sigma
S .
\end{equation}
Let $\alpha:\Sigma^qM\rightarrow M$ be the Adams map and $V(1)$ be
its cofibre given by the cofibration
\begin{equation}\label{1.2}
\Sigma^qM\stackrel{\alpha}{\rightarrow}M\stackrel{i^{\prime}}{\rightarrow}V(1)
\stackrel{j^{\prime}}{\rightarrow}\Sigma^{q+1}M.
\end{equation}
Let $\beta: \Sigma^{(p+1)q}V(1)\rightarrow V(1)$ be the $v_2$-map.
\\

{\bf Definition  1.1}\quad {\it We define, for $t\geq 1$, the
$\beta$-family $\beta_t=jj^\prime\beta^t i^\prime i \in
\pi_{q(tp+t-1)-2}(S)$. Here the maps $i$,  $j$,  $i^\prime$,
$j^\prime$, and $\beta$ are given as above.}

We have the following known result.

{\bf Theorem 1.1}\cite[Theorem 2.12]{Miller}\quad {\it $\beta_t\neq0
\in \pi_\ast(S)$ for $p\geq 5$ and $t\geq 1$.}

To determine the stable homotopy groups of spheres is very
difficult. Thus not so many families of homotopy elements in the
stable homotopy groups of spheres have been detected. See, for
example, \cite{Lee, Liulevicius, RCohen}.

In \cite{Liu2}, X. Liu obtained the following theorem, which is
called the representative theorem.

{\bf Theorem 1.2}\cite[Theorem 1.3]{Liu2}\quad {\it\label{Theorem:a}
 For $ p\geq 5$ and $2\leq s<p$, there exists the second Greek  letter
element $$ \widetilde{\beta }_s \in {\rm Ext}_A^{s,
q(sp+s-1)+s-2}(\mathbb{Z}_p,  \mathbb{Z}_p),$$ which converges to
the $\beta$-family $\beta_s \in \pi_{q(sp+s-1)-2}(S)$ in the ASS.
Moreover, $\widetilde{\beta}_s$ is represented by
$$a_{2}^{s-2}h_{2,0}h_{1,1}\in E_1^{s,q(sp+s-1)+s-2, \ast}$$ in the
May spectral sequence (MSS).}\\

In \cite{Lin}, J. Lin detected a new family of stable homotopy
groups of spheres and showed the following theorem.

{\bf Theorem 1.3}\cite{Lin}\quad {\it \label{Theorem:h} For $p \geq
5$, $n \geq m+2 \geq 4$. Then  $$h_0 h_n h_m\in {\rm Ext}_A^{3,
q(p^n+p^m+1) }(\mathbb{Z}_p, \mathbb{Z}_p) $$ is a permanent cycle
in the ASS and it converges to a  family  of homotopy elements of
order $p$, denoted by $\xi_{m,n}$, in  the stable homotopy groups of
spheres
 $\pi_{q(p^n+p^m+1)-3}(S)$.}

In this paper, we consider the non-triviality of the composite
$\xi_{m,n}\beta_s$ and obtain the following theorem.

{\bf Theorem 1.4}\quad {\it  \label{Theorem:main} Let $p \geq 5$, $n
\geq m+2> 5$, $2 \leq s <p$. Then the product
$$h_0h_nh_m\widetilde{\beta }_s \neq 0 \in {\rm Ext}_A^{s+3,  t(s)+s-2}(\mathbb{Z}_p ,  \mathbb{Z}_p)$$ is a permanent
cycle in the ASS and converges to a nontrivial family of homotopy
elements $\xi_{m,n}\beta_s \in \pi_{t(s)+s-5}(S)$, where
$t(s)=q(p^n+p^m+sp+s)$.}

In this paper we make use of the ASS and the MSS to prove our
theorem, especially the MSS. The method of the proof is very
elementary. By this method, one can consider some similar problems,
for example, the non-triviality of the composite
$\xi_{m,n}\gamma_s$, where $\gamma_s$ is the known $\gamma$-family
(cf. \cite{Liu3}).

The paper is arranged as follows: after giving some important lemmas
on the MSS in Section 2, we will prove Theorem 1.4 in Section 3.

\section{The ASS and some lemmas on the MSS}
 One of the main tools to determine the stable homotopy
groups of spheres $\pi _\ast (S)$ is the ASS. In 1957, Adams
constructed such a machinery in the form of a spectral sequence that
making the doubly graded group ${\rm Ext}_A^{\ast ,
\ast}(\mathbb{Z}_p, \mathbb{Z}_p)$ to the $p$-primary components of
the stable homotopy groups of spheres by adapting the methods of
homological algebra. From then on, the ASS has been a powerful tool
in studying stable homotopy theory.

Let $X$ a spectrum of finite type and $Y$ a finite dimensional
spectrum. Then there is a natural spectral sequence $\{E_r^{s, t},
d_r\}$ which is called Adams spectral sequence and
\begin{equation}
E_2^{s,  t}={\rm Ext}_A^{s, t}(H^\ast( X;  \mathbb{Z}_p  ),  H^\ast
(Y;  \mathbb{Z}_p )) \Rightarrow ([Y, X]_{t - s} )_p,
\end{equation}
 where the differential is
 \begin{equation}
 d_r:  E_r^{s,  t}\to E_r^{s+r,  t+r-1}.
 \end{equation}
If $X$ and $Y$ are sphere spectra $S$, then in the ASS
\begin{equation}
E_2^{s,  t}={\rm Ext}_A^{s,  t}(\mathbb{Z}_p ,
\mathbb{Z}_p)\Rightarrow(\pi_{t-s}(S))_p ,
\end{equation}
 the
$p$-primary components of the group $\pi_{t-s}(S).$

There are three problems in using the ASS: the calculation of the
$E_2$-term, the computation of the differentials and the
determination of the nontrivial extensions from $E_{\infty}$ to
$\pi_{\ast}(S)$. So, in order to compute the stable homotopy groups
of spheres with the ASS, we must compute the $E_2$-term of the ASS,
${\rm Ext}_A^{\ast, \ast}(\mathbb{Z}_p, \mathbb{Z}_p)$. The most
successful tool for computing ${\rm Ext}_A^{\ast,\ast}(\mathbb{Z}_p,
\mathbb{Z}_p)$ is the MSS.

From \cite{Ravenel}, there is a MSS
 $\{E_r^{s, t, \ast}, d_r\}$
which converges to ${\rm Ext}_A^{s, t}(\mathbb{Z}_p, \mathbb{Z}_p)$
with $E_1$-term
\begin{equation}\label{2.4}
E_1^{\ast, \ast, \ast}= E(h_{m, i}|m>0, i\geq 0)\otimes P(b_{m,
i}|m>0, i\geq 0)\otimes P(a_n|n\geq 0),
\end{equation}
 where $E$ is the exterior algebra, $P$ is the polynomial
algebra, and
$$h_{m,  i}\in E_1^{1,  2(p^m-1)p^i,  2m-1},   b_{m,  i}\in
E_1^{2,  2(p^m-1)p^{i+1},  p(2m-1)},  a_n\in E_1^{1,  2p^n-1,
2n+1}.$$ The $r$-th May differential is
\begin{equation}\label{2.5}
d_r:  E_r^{s,  t, u}\rightarrow E_r^{s+1,  t,  u-r},
 \end{equation}
and if $x\in E_r^{s,  t,  \ast}$ and $y\in E_r^{s^{\prime},
t^{\prime},  \ast}$, then $d_r(x\cdot y)=d_r(x)\cdot y+(-1)^sx\cdot
d_r(y)$. From \cite[Proposition 2.5]{Liuwang}, there exists a graded
commutativity in the May $E_1$-term as follows:
\begin{equation}\label{2.6}\left\{\begin{array}{ll}
a_mh_{n,j}=h_{n,j}a_m, & h_{m,k}h_{n,j}=-h_{n,j}h_{m,k},\\
  a_mb_{n,j}=b_{n,j}a_m,&
h_{m,k}b_{n,j}=b_{n,j}h_{m,k},\\ a_ma_n=a_na_m, & b_{m,n}b_{i,j}=
b_{i,j}b_{m,n}.
\end{array}\right.
\end{equation} The first May
differential $d_1$ is given by
\begin{equation}\label{2.7}
\left\{\begin{array}{l}
d_1(h_{i,  j})=\sum\limits_{0<k<i}^{}h_{i-k,  k+j}h_{k, j},\\
d_1(a_i)=\sum\limits_{0\leq k<i}h_{i-k,  k}a_k,\\
d_1(b_{i,  j})=0.
\end{array}
\right.
\end{equation}

For each element $x\in E_{1}^{s, t, \mu}$, we define
$\hbox{filt}~x=s$, $\hbox{deg}~x=t$, $\hbox{M}~(x)=\mu$. Then we
have
\begin{equation}\label{2.8}
\left\{\begin{array}{l}
\hbox{filt}~h_{i, j}=\hbox{filt}~a_{i}=1,\   \hbox{filt}~b_{i, j}=2,\\
\hbox{deg}~h_{i, j}=2(p^i-1)p^j=q(p^{i+j-1}+\cdots+p^j),\\
\hbox{deg}~b_{i, j}=2(p^i-1)p^{j+1}=q(p^{i+j}+\cdots+p^{j+1}),\\
\hbox{deg}~a_i=2p^i-1=q(p^{i-1}+\cdots+1)+1,\\
\hbox{deg}~a_0=1, \\
\hbox{M}~(h_{i, j})=\hbox{M}~(a_{i-1})=2i-1,\\
\hbox{M}~(b_{i, j})=(2i-1)p, \end{array}\right. \end{equation}
 where $i\geq 1$, $j\geq 0$.

In Section 3, we will need the following lemmas on the MSS.

By the knowledge on $p$-adic expression in number theory, we have
that for each integer $t\geq 0$, it can be always expressed uniquely
as  $$t=q({c_n p^n+c_{n-1} p^{n-1}+\cdots +c_1
p+c_0})+c_{-1},$$where $0\leq c_i <p$ $(0\leq i <n)$, $0< c_n
<p$,  $0 \leq c_{-1}<q$.\\

{\bf Lemma  2.1}\cite[Proposition 1.1]{Liu3}\quad {\it
\label{Lemma:j} Let $t$ as above. Let $s_1$ be a positive integer
with $0< s_1<p$. If there exists some $0\leq j \leq n$ such that
$c_j>s_1$, then in the MSS $$ E_1^{s_1, t, \ast}=0.$$ } \hfill$\Box$

{\bf Lemma  2.2}\quad {\it \label{Lemma:jjj} Let $t$ as above. Let
$s_1$ be a positive integer with $0<s_1<q$. If $c_{-1}>s_1$, then in
the MSS,  $$E_1^{s_1,t,\ast}=0.$$ }

{\bf Proof}\quad The proof is similar to that of \cite[Proposition
1.1]{Liu3} and is omitted here. \hfill$\Box$\\

Let $t$ as above and $s$ a given positive integer. Suppose that in
the MSS a generator $\omega \in E_1^{s,t,\ast}$ is of the form
$w=x_1x_2\cdots x_m$, where $x_i$ is one of $a_k$, $h_{l,j}$ or
$b_{u,z}$, $1 \leq i \leq m$, $0\leq k \leq
 n+1$, $0<u+z\leq n$, $0<l+j\leq n+1$, $l>0$, $j\geq 0$,
 $u>0$, $z\geq 0$. By (2.8), we
can assume that for any $1\leq i\leq m$ $\hbox{deg}~ x_i=q(c_{i,
n}p^n+c_{i, n-1}p^{n-1}+\cdots+c_{i,1}p+c_{i, 0})+{c}_{i,-1} $,
where $c_{i, j}=0$ or $1$ for $0 \leq j \leq n$, ${c}_{i,-1}=1$ if
$x_i=a_{k_{i}}$, or ${c}_{i,-1}=0$. It follows that $$\hbox{deg}~
\omega=\sum\limits_{i=1}^{m}{\hbox{deg}~
x_i}=q[(\sum\limits_{i=1}^{m}{c_{i, n}})p^n+\cdots
+(\sum\limits_{i=1}^{m}{c_{i, 1}})p^{1}+\sum\limits_{i=1}^{m}{c_{i,
0}}]+\sum\limits_{i=1}^{m}{c_{i, -1}}.$$ For convenience, we denote
$\sum\limits_{i=1}^{m}{c_{i, j}}$ by $\bar{c}_j$ for $j\geq -1$.

{\bf Lemma 2.3}\quad {\it With notation as above. If there exist
three integers $-1\leq i_1<i_2<i_3\leq n$ such that
$\bar{c}_{i_1}+\bar{c}_{i_3}-m>\bar{c}_{i_2}$, then $w$ is
impossible to exist.}

{\bf Proof}\quad  By (\ref{2.8}) and (\ref{2.4}), one easily gets
the
lemma.\hfill$\Box$\\

{\bf Lemma  2.4}\quad {\it With notation as above. Suppose that
$m=s$, and there exist three integers $i_1$, $i_2$ and $i_3$
satisfying the following conditions that

 {\rm (i)} $-1 \leq i_1  <i_2<i_3\leq n${\rm ;}

  {\rm (ii)} $\overline{c}_{i_1}+\overline{c}_{i_3}-m \leq
\overline{c}_{i_2}${\rm ;}

{\rm (iii)} $\overline{c}_{j}= \left\{\begin{array}{cc}
0& {-1\leq j<i_1}\\
0& {i_3<j\leq n}.
\end{array}\right. $ \\
Then we have the following consequences:

{\rm (1)} When $i_1>-1$, there are
$(\overline{c}_{i_{1}}+\overline{c}_{i_{3}}-m)$
$h_{{\overline{c}_{i_{3}}-\overline{c}_{i_{1}}+1},\overline{c}_{i_{1}}}$'s
among $\omega$. Furthermore, if
$\overline{c}_{i_{1}}+\overline{c}_{i_{3}}-m>1$, then $w=0$.

{\rm (2)} When $i_1=-1$, there are
$(\overline{c}_{i_{1}}+\overline{c}_{i_{3}}-m)$ $a_{i_{3}+1}$'s
among $\omega$.}

{\bf Proof }\quad By (\ref{2.8}) and (\ref{2.4}), the desired
results easily follow.\hfill$\Box$

\section{Proof of Theorem 1.4} In this section, we will determine
two ${\rm Ext}$ groups which will be used in the proof of Theorem
1.4. In order to do it, we first consider some May $E_1$-terms $E_1^{u,v,\ast}$ with two given integers $u$ and $v$, and show the following lemma.\\

{\bf Lemma 3.1} \quad {\it \label{Lemma:d}
  Let $p \geq 5$, $n \geq m+2> 5$, $2 \leq s <p$ and  $1\leq r \leq s+3$. Then in the
  MSS, we have
\begin{equation}\label{3.1}
E_{1}^{{s+3-r}, {t(s)+s-r-1}, \ast}= \left\{\begin{array}{ll}
\mathbb{Z}_p\{\mathbf{g}_1,\cdots,\mathbf{g}_7\}& {r=1 \ and \ s=p-1,}\\
0& other.
\end{array}\right.
\end{equation}
Here, $t(s)=q(p^n+p^m+sp+s)$, and $\mathbf{g}_1$, $\cdots$,
$\mathbf{g}_7$ equal elements $a_n^{p-3}h_{3, 0}h_{1, m}h_{n-2,
2}h_{n, 0}$, $a_n^{p-3}h_{1,2}h_{m+1,0}h_{n-m, m}h_{n, 0}$,~~
$a_{m+1}a_n^{p-4}h_{3, 0}h_{n-m, m}h_{n-2, 2}h_{n, 0}$,~~
$a_n^{p-3}h_{3, 0}h_{m-1, 2}h_{n-m, m}h_{n, 0}$,\\
$a_n^{p-3}h_{3, 0}h_{m+1, 0}h_{n-m, m}h_{n-2, 2}$,
 $a_{3}a_n^{p-4}  h_{m+1, 0}h_{n-m,
m}h_{n-2, 2}h_{n, 0}$ and $a_m^{p-3}h_{3, 0}h_{m, 0}h_{m-2, 2}h_{1,
n}$, respectively.}

 {\bf Proof} \quad We divide the proof into the following two cases.

 {\bf Case 1} $s-r-1<0$. By the knowledge on $p$-adic
 expression in number theory and $1 \leq r \leq s+3$,
 we would have $s+3-r <s-r-1+q<q$. In this case
 $$E_{1}^{{s+3-r}, {t(s)+s-r-1}, \ast}=0$$ by Lemma $2.2$.

{\bf Case 2} $s-r-1\geq 0$. Thus $1\leq r\leq s-1$. If $r\geq 4$,
then $s+3-r<s$, which implies that in this case $E_{1}^{{s+3-r},
{t(s)+s-r-1}, \ast}=0$ by Lemma 2.1. Consequently, in the rest of
the proof, we always assume {$r\leq
 3$}.

 Consider $\omega=x_1 x_2 \cdots x_{m^\prime}\in E_{1}^{s+3-r,
t(s)+s-r-1, \ast}$
 in the MSS, where $x_i $ is one of $ a_k$, $h_{l,j}$, $b_{u,z}$, $1
\leq i \leq m^{\prime}$, $0\leq k \leq
 n+1$,
 $0<l+j\leq
n+1$, $0<u+z\leq n$, $l>0$, $j \geq 0$, $u>0$, $z \geq 0$. By
$(\ref{2.8})$, we can assume that $\hbox{deg}~ x_i=q(c_{i,
n}p^n+c_{i, n-1}p^{n-1}+\cdots+c_{i,1}p+c_{i, 0})+{c}_{i,-1} $,
where $c_{i, j}=0$ or $1$ for $0\leq j \leq n$, ${c}_{i,-1}=1$ if
$x_i=a_{k_{i}}$, or ${c}_{i,-1}=0$. It follows that
\begin{equation} \label{equation:b}
\left\{\begin{array}{l}
 \hbox{filt}~ \omega =\sum\limits_{i=1}^{m^{\prime}}\hbox{filt}~ { x_i} =s+3-r,\\
\hbox{deg}~  \omega=\sum\limits_{i=1}^{m^\prime}{\hbox{deg}~
x_i}=q[(\sum\limits_{i=1}^{m^\prime}{c_{i,
n}})p^n+(\sum\limits_{i=1}^{m^\prime}{c_{i, n-1}})p^{n-1}+\cdots
+(\sum\limits_{i=1}^{m^\prime}{c_{i,  m}})p^{m}\\

{\phantom {\hbox{deg}~ \omega +}+(\sum\limits_{i=1}^{m^\prime}{c_{i,
m-1}})p^{m-1}+\cdots
 +(\sum\limits_{i=1}^{m^\prime}{c_{i, 1}})p+(\sum\limits_{i=1}^{m^\prime}{c_{i, 0}})]+
 (\sum\limits_{i=1}^{m^\prime}{c_{i,-1}})}\\
{\phantom {\hbox{deg}~ \omega} =t(s)+s-r-1.}
\end{array}\right.
\end{equation}
Note that $\hbox{filt}~x_i=1$ or $2$ and $2\leq s<p$. From
$\sum\limits_{i=1}^{m^{\prime}} \hbox{filt}~x_i=s+3-r$, it follows
that $$m^{\prime}\leq s+2<p+2.$$
 Using $0 \leq {s,\ s-r-1}<p$ and  the knowledge on the $p$-adic
 expression in number theory,  we have the following equations from
 (\ref{equation:b}).

{\begin{equation}\label{equation:d}
 \left\{\begin{array}{ll}
\sum\limits_{i=1}^{m^\prime}{c}_{i,-1}=s-r-1+\lambda_{-1}q,&
 \lambda_{-1} \geq 0,\\
\sum\limits_{i=1}^{m^\prime}c_{i,  0}+\lambda_{-1}=s+\lambda_{0}p,
 &\lambda_{0} \geq 0,\\
\sum\limits_{i=1}^{m^\prime}c_{i,  1}+\lambda_{0}=s+\lambda_{1}p,
 &\lambda_{1} \geq  0,\\
\sum\limits_{i=1}^{m^\prime}c_{i,  2}+\lambda_{1}=0+\lambda_{2}p,
 &\lambda_{2} \geq 0,\\
\sum\limits_{i=1}^{m^\prime}c_{i,  3}+\lambda_{2}=0+\lambda_{3}p,
 &\lambda_{3} \geq 0,\\
 \cdots &\cdots\\
\sum\limits_{i=1}^{m^\prime}c_{i,
m-1}+\lambda_{m-2}=0+\lambda_{m-1}p,&
 \lambda_{m-1} \geq 0,\\
 \sum\limits_{i=1}^{m^\prime}c_{i,  m}+\lambda_{m-1}=1+\lambda_{m}p,
& \lambda_{m} \geq  0,\\
\sum\limits_{i=1}^{m^\prime}c_{i, m+1}+\lambda_{m}=0+\lambda_{m+1}p,
& \lambda_{m+1} \geq 0,\\
 \cdots & \cdots\\
\sum\limits_{i=1}^{m^\prime}c_{i,
n-2}+\lambda_{n-3}=0+\lambda_{n-2}p,&
 \lambda_{n-2} \geq 0,\\
 \sum\limits_{i=1}^{m^\prime}c_{i,
n-1}+\lambda_{n-2}=0+\lambda_{n-1}p,&
 \lambda_{n-1} \geq 0,\\
\sum\limits_{i=1}^{m^\prime}c_{i,  n}+\lambda_{n-1}=1.&

\end{array}\right.
\end{equation}
} By the knowledge on the $p$-adic expression and $m^\prime <p+2$,
we
 have $\lambda_{-1}=\lambda_{0}=\lambda_{1}=0$. For
convenience, in the rest of the proof we will use $\overline{c}_{j}$
to denote $\sum\limits_{i=1}^{m^\prime}{c_{i, j}} $  for $-1 \leq j
\leq n$. From  the fourth equation of (\ref{equation:d})
$\overline{c}_2=\lambda_2 p$, $\lambda_2$ may equal $0$ or $1$.

 {\bf Subcase $2.1$ } $\lambda_2=0$.

{\bf Assertion 3.1}\quad { If $\lambda_{2}=0$, then $\lambda_{3}=
\cdots =\lambda_{m-1}=0$. }

Suppose $\lambda_3=1$. Then from the fifth equation of
(\ref{equation:d}) we would have $\overline{c}_3=p$, which implies
that $m^{\prime}$ can only equal $p$ or $p+1$. Note that $2\leq s
<p$. From $\overline{c}_{3}=p$, $\overline{c}_{2}=0$ and
$\overline{c}_{1}=s$, one would have
$\overline{c}_{3}+\overline{c}_{1}-m^{\prime}=p+s-m^{\prime} \geq
1>0=\overline{c}_{2}$. Thus by Lemma $2.3$, $\omega$ is impossible
to exist. Thus, $\lambda_3=0$. Similarly, one can show that
$\lambda_{4}= \cdots =\lambda_{m-1}=0$. Assertion 3.1 is proved.

From the $(m+2)$-th equation of ($3.3$) $\overline{c}_m=1+\lambda_m
p$, $\lambda_m$ may equal $0$ or $1$.

{\bf  Subcase $2.1.1$ } $ \lambda_m =0$. An argument similar to that
used in Assertion $3.1$ shows that $\lambda_{m+1}=\cdots
=\lambda_{n-1}=0$. Thus we have
\begin{center}
\begin{tabular}{|c|c|c|c|c|c|c|c|c|c|c|c|}
\hline
$\overline{c}_{n}$& $\overline{c}_{n-1}$ & $\cdots$ & $\overline{c}_{m+1}$&$\overline{c}_{m}$&$\overline{c}_{m-1}$&$\cdots $&$\overline{c}_{3}$&$\overline{c}_{2}$&$\overline{c}_{1}$&$\overline{c}_{0}$&$\overline{c}_{-1}$\\
\hline $1$& $0$& $\cdots$& $0$& $1$&$0$ &$\cdots $ &$0$ &$0$ &$s$ & $s$ &$s-r-1$ \\
 \hline
\end{tabular}.
\end{center}

If $\omega$ has $h_{1,n}h_{1,m}$ as factors, one can let
$\omega=h_{1,n}h_{1,m}\omega_1$ by (\ref{2.6}). Then $\hbox{filt}~
\omega_1=s+1-r$ and $\hbox{deg}~ \omega_1= spq+sq+(s-r-1)$. When $r
> 1$, $\omega_1$ is impossible to exist by
Lemma $2.1$. So $\omega$ is impossible to exist either. When $r=1$,
$\omega_1$ has $(s-2)$ $a_2$'s among $\omega$ if $\omega$ exists by
Lemma $2.4$. Then up to sign $\omega_1=a_2^{s-2}\omega_2$ with
$\omega_2 \in E_{1}^{2,2pq+2q,\ast}=0$, which  means $\omega=0$.

Similarly, $\omega$ cannot have $h_{1,n}b_{1,m-1}$,
$b_{1,n-1}h_{1,m}$, $b_{1,n-1}b_{1,m-1}$ as factors either.

{\bf Subcase 2.1.2}\quad $\lambda_{m}=1$. In this case
$\lambda_{m+1}=\cdots=\lambda_{n-1}=1$. An argument similar to that
used in Assertion 3.1 can show that in this case $\omega$ is
impossible to exist.

{\bf  Subcase $2.2$ } $ \lambda_2 =1$. In this case,
$(\lambda_{3},\cdots,\lambda_{m-1})$ must equal $(1,\cdots,1)$. From
the $(m+2)$-th equation of (\ref{equation:d})
$\overline{c}_m=\lambda_m p$, $\lambda_m$ may equal $0$ or $1$.

{\bf  Subcase $ 2.2.1$}\quad $\lambda_{m}=1$. In this case,
$(\lambda_{m+1},\cdots,\lambda_{n-1})$ must equal $(1,\cdots,1)$.
Thus we have
\begin{center}
\begin{tabular}{|c|c|c|c|c|c|c|c|c|c|c|c|}
\hline
$\overline{c}_{n}$& $\overline{c}_{n-1}$ & $\cdots$ & $\overline{c}_{m+1}$&$\overline{c}_{m}$&$\overline{c}_{m-1}$&$\cdots $&$\overline{c}_{3}$&$\overline{c}_{2}$&$\overline{c}_{1}$&$\overline{c}_{0}$&$\overline{c}_{-1}$\\
\hline $0$& $p-1$& $\cdots$& $p-1$& $p$&$p-1$ &$\cdots $ &$p-1$ &$p$ &$s$ & $s$ &$s-r-1$ \\
\hline
\end{tabular}.
\end{center}

If $r=2$ or $3$, then by $p\geq 5$, $2\leq s<p$ and $m^{\prime}\leq
s+3-r$ one can have
$$\overline{c}_{2}+\overline{c}_{m}-m'=p+p-m^{\prime}\geq
p+p-(s+1)\geq p>p-1=\overline{c}_{3},$$ which implies that $\omega$
is impossible to exist by Lemma $2.3$.

If $r=1$, then one has $\hbox{filt}~\omega=s+2$. From
$\overline{c}_{m}=p$, one has $m'\geq p$ by $c_{i,m}=0$ or $1$. Thus
$m^{\prime}$ may equal $p$ or $p+1$. If $m^{\prime}=p$, then
$\overline{c}_2+\overline{c}_m-m^{\prime}=p>p-1=\overline{c}_{3}$,
which implies that $\omega$ is impossible to exist by Lemma 2.3.
Thus, in the rest of Subcase $2.2.1$ we always assume that {$r=1$
and $m^{\prime}=p+1$}. Thus we have {$s=p-1$},
$\hbox{filt}~\omega=p+1$ and $\omega=x_1 \ \cdots \ x_{p+1}\in
E(h_{i,j}|i>0,j\geq 0) \otimes P(a_n|n \geq 0)$. The table above
becomes
\begin{center}
\begin{tabular}{|c|c|c|c|c|c|c|c|c|c|c|c|}
\hline
$\overline{c}_{n}$& $\overline{c}_{n-1}$ & $\cdots$ & $\overline{c}_{m+1}$&$\overline{c}_{m}$&$\overline{c}_{m-1}$&$\cdots $&$\overline{c}_{3}$&$\overline{c}_{2}$&$\overline{c}_{1}$&$\overline{c}_{0}$&$\overline{c}_{-1}$\\
\hline $0$& $p-1$& $\cdots$& $p-1$& $p$&$p-1$ &$\cdots $ &$p-1$ &$p$ &$p-1$ & $p-1$ &$p-3$ \\
 \hline
\end{tabular}.
\end{center}

{\bf Assertion 3.2}\quad {$\omega$ has $p-1$ factors whose degrees
are $q($higher terms  on $ p+p^m+\cdots +p^2+$ lower terms on
$p)+\epsilon$, where $\epsilon =0$ or $1$, and two factors whose
degree are $q($ higher terms on $p+p^m)$ and $q(p^2 +$ lower terms
on $p )+\epsilon $, respectively.}

This assertion can be easily verified by (\ref{2.8}) and
(\ref{equation:b}).

{\bf Assertion 3.3 }\quad {$\omega $ cannot have $h_{2,1}$ or
$h_{j,m}$ ($2 \leq j <n-m$) as a factor. }

Otherwise, we can let $\omega=\omega_1h_{2,1}$ by (\ref{2.6}). Then
$\hbox{filt}~\omega_1=p$, $\hbox{deg}~\omega_1=q[(p-1)p^{n-1}+\cdots
+(p-1)p^{m+1}+pp^{m}+(p-1)p^{m-1}+\cdots
+(p-1)p^{3}+(p-1)p^{2}+(p-2)p+(p-1)]+p-3$. In this case $\omega_1$
is impossible to exist by Lemma $2.3$. Thus $\omega$ cannot have
$h_{2,1}$ as a factor. Similarly, $\omega$ cannot have $h_{j,m}$ $(2
\leq j <n-m)$ as a factor.

From Assertions 3.2 and 3.3, there must be one of $h_{1,2}h_{1,m}$,
$h_{3,0}h_{1,m}$, $a_3 h_{1,m}$, $h_{1,2}h_{n-m,m}$,
$h_{3,0}h_{n-m,m}$ or $a_3 h_{n-m,m}$ among $\omega$ if $\omega$
exists. By (\ref{2.6}), we let $\omega=\omega_1 \omega_2$, where
$\omega_2$ is one of the six factors above. Then
$\hbox{filt}~\omega_1=p-1$.

{\bf   (i)}\quad If $\omega_2=h_{1,2}h_{1,m}$, then
$\hbox{deg}~\omega_1=q
[(p-1)p^{n-1}+\cdots+(p-1)p^m+\cdots+(p-1)p+(p-1)]+p-3$. So there
must be two $h_{n,0}$'s in $\omega$ by  Lemma 2.4 (1), which implies
that $\omega_1=0$. Then $\omega=0$.

Similarly, one can show that $\omega=0$ if $\omega_2=a_3 h_{1,m}$.

{\bf   (ii)}\quad If $\omega_2=h_{3,0}h_{1,m}$, then
$\hbox{deg}~\omega_1=q
[(p-1)p^{n-1}+\cdots+(p-1)p^m+\cdots+(p-1)p^2+(p-2)p+(p-2)]+p-3$. By
Lemma 2.4, $\omega_1$ must equal $a_n^{p-3}h_{n,0}h_{n-2,2}$ up to
sign . Thus up to sign $\omega=a_n^{p-3}h_{3, 0}h_{1, m}h_{n-2,
2}h_{n, 0}$, denoted by $\mathbf{g}_1$.

{\bf   (iii)}\quad If $\omega_2=h_{1,2}h_{n-m,m}$, then
$\hbox{deg}~\omega_1=q
[(p-2)p^{n-1}+\cdots+(p-2)p^{m+1}+(p-1)p^m+(p-1)p^{m-1}+\cdots+(p-1)p^2+(p-1)p+(p-1)]+p-3$,
so $\omega_1$ has at least $p-4$ $a_n$'s by Lemma $2.4$. We let
$\omega_1=\omega_3 a_n^{p-4}$  by (\ref{2.6}). Thus
$\hbox{filt}~\omega_3=3$,
$\hbox{deg}~\omega_3=q(2p^{n-1}+\cdots+2p^{m+1}+3p^{m}+3p^{m-1}+\cdots+3p^3+3p^2+3p+3)+1$.
Then $\omega_3\in E_{1}^{3,\hbox{deg}~\omega_3,\ast}=\mathbb{Z}_p
\{a_n  h_{n,0} h_{m+1,0} \}$. Thus up to sign
$\omega=a_n^{p-3}h_{1,2}h_{m+1,0}h_{n-m, m}h_{n, 0}$, denoted by
$\mathbf{g}_2$.

{\bf   (iv)}\quad If $\omega_2=h_{3,0}h_{n-m,m}$, an argument
similar to that used in (iii) shows $\omega_1= a_n^{p-4}\omega_3$
with $\omega_3 \in E_{1}^{3,t,\ast}=\mathbb{Z}_p \{
a_{m+1}h_{n-2,2}h_{n,0}, a_n h_{m-1,2} h_{n,0}, a_n  h_{m+1,0}
h_{n-2,2} \}$, where
$t=q(2p^{n-1}+\cdots+2p^{m+1}+3p^{m}+\cdots+3p^{2}+2p+2)+1$. Thus up
to sign $\omega=a_{m+1}a_n^{p-4} h_{3, 0}h_{n-m, m}h_{n-2, 2}h_{n,
0}$, $a_n^{p-3}h_{3, 0}h_{m-1, 2}h_{n-m, m}h_{n, 0}$ or
$a_n^{p-3}h_{3, 0}h_{m+1, 0}h_{n-m, m}h_{n-2, 2} $, denoted by
$\mathbf{g}_3$, $\mathbf{g}_4$, $\mathbf{g}_5$ respectively.

{\bf (v)}\quad If $\omega_2=a_3 h_{n-m,m}$, by an argument similar
to that used in (iii) we have $\omega_1=a_n^{p-5} \omega_3$ with
$\omega_3\in E_{1}^{4,t',\ast}=\mathbb{Z}_p \{a_n h_{m+1,0} h_{n,0}
h_{n-2,2} \}$, where
$t'=q(3p^{n-1}+\cdots+3p^{m+1}+4p^{m}+\cdots+4p^{2}+3p+3)+1$. Thus
up to sign $\omega=a_{3}a_n^{p-4}  h_{m+1, 0}h_{n-m, m}h_{n-2,
2}h_{n, 0}$, denoted by $\mathbf{g}_6$.

{\bf Subcase  $ 2.2.2$}\quad  $ \lambda_m =0$. By an argument
similar to that used in Assertion $3.1$, we have $\lambda_{m+1}=
\cdots =\lambda_{n-1}=0$. Thus we have
\begin{center}
\begin{tabular}{|c|c|c|c|c|c|c|c|c|c|c|c|}
\hline $\overline{c}_{n}$& $\overline{c}_{n-1}$ & $\cdots$ &
$\overline{c}_{m+1}$&$\overline{c}_{m}$&$
\overline{c}_{m-1}$&$\cdots $&$\overline{c}_{3}$&$\overline{c}_{2}$&$\overline{c}_{1}$&$\overline{c}_{0}$&$\overline{c}_{-1}$\\
\hline $1$& $0$& $\cdots$& $0$& $0$&$p-1$ &$\cdots $ &$p-1$ &$p$ &$s$ & $s$ &$s-r-1$ \\
 \hline
\end{tabular}.
\end{center}
Obviously $\omega$ must have a factor $h_{1,n}$. One can let
$\omega=\omega_1 h_{1,n}$ by (\ref{2.6}).

If $r=2$ or $3$, it is easy to show that $\omega_1$ is impossible to
exist by Lemma 2.1, which implies that $\omega$ is impossible to
exist either.

If $r=1$, by Lemma 2.1 it is easy to get that in this case $s$ must
equal $p-1$ and $m^{\prime}$ must equal $p+1$. By an argument
similar to that used in  Subcase $2.2.1$, we get that up to sign
$\omega=a_m^{p-3}h_{3, 0}h_{m, 0}h_{m-2, 2}h_{1, n} \in E^{p+1,
t(p-1)+p-3, (2m+1)p-2m-3},$ denoted by $\mathbf{g}_7$.

Combining Cases 1 and 2, we complete the proof of the lemma. \hfill$\Box$\\

By use of Lemma 3.1, we now show the non-triviality of $h_0 h_n h_m
\widetilde{\beta }_s$.

{\bf Theorem 3.2}\quad {\it \label{Theorem:xx} Let $p \geq 5$, $n
\geq m+2> 5$, $2 \leq s <p$. Then the product
 $$ h_0 h_n h_m \widetilde{\beta }_s \neq 0 \in {\rm Ext}_{A}^{{s+3}, {t(s)+s-2}}(\mathbb{Z}_p ,
 \mathbb{Z}_p),$$ where $t(s)=q(p^n+p^m+sp+s)$}.

{\bf Proof  }\quad Since $h_{1, n}(n \geq 0)$ and $a_{2}^{s-2}h_{2,
0} h_{1, 1}$ are permanent cycles in the MSS and converge
nontrivially to $h_n$,
 $\widetilde{\beta_s} \in {\rm Ext}_A^{\ast, \ast}(\mathbb{Z}_p,
\mathbb{Z}_p)$ respectively, $a_{2}^{s-2}h_{2, 0}h_{1, 1}h_{1,
0}h_{1, n}h_{1, m} \in E_{1}^{{s+3}, {t(s)+s-2}, \ast} $ is a
permanent cycle in the MSS and converges to $h_0 h_n h_m
\widetilde{\beta }_s \in {\rm Ext}_{A}^{{s+3},
{t(s)+s-2}}(\mathbb{Z}_p, \mathbb{Z}_p).$

{\bf Case $1$} $s=p-1$. By (\ref{2.7}), one can have that up to sign
\begin{equation}\label{3.4}
\left\{ \begin{array}{ll}
d_1(\mathbf{g}_1)=a_n^{p-3}h_{1, 0}h_{3, 0}h_{1, m}h_{n-2, 2}h_{n-1, 1}+\cdots &\neq 0;\\
d_1(\mathbf{g}_2)=a_n^{p-3}h_{1, 0}h_{1, 2}h_{m+1, 1}h_{n-m, m}h_{n-1, 1}+\cdots &\neq 0;\\
d_1(\mathbf{g}_3)=a_n^{p-4}a_{m+1}h_{1, 0}h_{3, 0}h_{n-m, m}h_{n-2, 2}h_{n-1, 1}+\cdots &\neq 0; \\
d_1(\mathbf{g}_4)=a_n^{p-3}h_{1, 0}h_{3, 0}h_{m-1, 2}h_{n-m, m}h_{n-1, 1}+\cdots  &\neq 0; \\
d_1(\mathbf{g}_5)=a_n^{p-3}h_{1, 0}h_{3, 0}h_{m, 1}h_{n-m, m}h_{n-2, 2}+\cdots &\neq 0;\\
d_1(\mathbf{g}_6)=a_n^{p-4}a_3 h_{1, 0}  h_{m+1, 0}h_{n-m, m}h_{n-2, 2}h_{n-1, 1}+\cdots &\neq 0;\\
d_1(\mathbf{g}_7)= a_m^{p-3}h_{1, 2}h_{3, 0}h_{m-3, 3}h_{1, n}h_{n,
0}+\cdots &\neq 0.
\end{array} \right. \end{equation}
Obviously the first May differential of each of the seven generators
contains at least a term which is not in the first May differentials
of the other generators, which implies that $d_1(\mathbf{g}_1),\
\cdots,\ d_1(\mathbf{g}_7)$ are linearly independent. Thus,
$E_{2}^{{p+1}, {t(p-1)+p-3}, \ast}=0.$ It follows that
$$E_{r}^{{p+1}, {t(p-1)+p-3}, \ast}=0~{\rm for}~r\geq 2.$$
Meanwhile, by (\ref{2.8}) we have that
$M(\mathbf{g}_i)=(2n+1)p-2n-3$( $1\leq i \leq 6$),
$M(\mathbf{g}_7)=(2m+1)p-2m-3$ and
$M(a_2^{p-3}h_{2,0}h_{1,1}h_{1,0}h_{1,n}h_{1,m})=5p-8$. Then from
(\ref{2.5}) one has $a_2^{p-3}h_{2,0}h_{1,1}h_{1,0}h_{1,n}h_{1,m}
\notin d_1(E_1^{p+1, t(p-1)+p-3, p(2n+1)-2n-3})$ and
$a_2^{p-3}h_{2,0}h_{1,1}h_{1,0}h_{1,n}h_{1,m} \notin d_1(E_1^{p+1,
t(p-1)+p-3, p(2m+1)-2m-3})$. Thus we have the permanent cycle
$a_2^{p-3}h_{2,0}h_{1,1}h_{1,0}h_{1,n}h_{1,m} \in E_r^{p+2,
t(p-1)+p-3,\ast }$ cannot be hit by any May differential. It follows
that $h_0 h_n h_m \widetilde{\beta }_{p-1} \neq 0$.

{\bf Case $2$ }$2\leq s<p-1$. From Lemma 3.1 one has that in this
case the May $E_1$-term
 $$E_{1}^{{s+2}, {t(s)+s-2}, \ast}=0.$$ Thus one has
$$E_{r}^{{s+2}, {t(s)+s-2}, \ast}=0 \  {\rm for} \ r>1.$$ Consequently, the
permanent cycle $h_{1, 0}h_{1, n}h_{1, m}a_{2}^{s-2}h_{2, 0}h_{1, 1}
\in E_{r}^{{s+3}, {t(s)+s-2}, \ast} $ cannot be hit by any
differential in the MSS. Then $h_0 h_n h_m \widetilde {\beta_{s}}
\neq 0 \in {\rm Ext}_{A}^{{s+3}, t(s)+s-2}(\mathbb{Z}_p,
\mathbb{Z}_p).$

From Cases $1$ and $2$, we complete
 the proof of the theorem.\hfill$\Box$\\

{\bf Theorem 3.3}\quad {\it \label{Theorem:f} Let $p \geq 5$, $n
\geq m+2> 5$, $2 \leq s <p$ and $2\leq r \leq s+3$. Then
$${\rm Ext}_{A}^{{s+3-r},{t(s)+s-r-1}, \ast}(\mathbb{Z}_p,
\mathbb{Z}_p)=0,$$ where $t(s)=q(p^n+p^m+sp+s)$. }

{\bf Proof} \quad From the case $2\leq r \leq  s+3 $ in Lemma
 3.1, we have that in the MSS
$$E_{1}^{{s+3-r},{t(s)+s-r-1}, \ast}=0.$$
 The theorem
follows easily by the MSS.  \hfill$\Box$\\

Now we give the proof of Theorem $1.4$ .

{\bf Proof of Theorem 1.4}  \quad  From Theorem 1.2, we have that
$\widetilde{\beta}_s$ converges to $\beta$-family $\beta_s \in
\pi_{spq+(s-1)q+s-2}(S)$ in ASS. From Theorem $1.3$,  $h_0 h_n
h_m\in {\rm Ext}_A^{3, p^n q +p^m q+q }(\mathbb{Z}_p, \mathbb{Z}_p)$
is a permanent cycle in the ASS and converges to a nontrivial family
of homotopy elements $\xi_{m,n} \in \pi_{p^n q+p^m q-3}(S)$. Hence,
we have that the composite $$\xi_{m,n}\beta_{s}$$ is represented up
to a nonzero scalar by $$h_0 h_n h_m \widetilde{\beta }_s \neq 0 \in
{\rm Ext}_A^{s+3, t(s)+s-2}(\mathbb{Z}_p,  \mathbb{Z}_p)$$ in the
ASS (cf. Theorem $3.2$).

Moreover, from Theorem $3.3$,  $h_0 h_n h_m \widetilde{\beta }_s$
cannot be hit by any differential in the ASS. Consequently, the
corresponding homotopy element $\xi_{m,n}\beta_{s}$ is nontrivial.
This proves Theorem $1.4$. \hfill$\Box$

\bibliographystyle{amsplain}

\end{document}